\documentclass[12pt,a4paper,twoside]{amsart}
\usepackage{a4}

\usepackage[cp1251]{inputenc}
\usepackage[T2A]{fontenc}

\usepackage[russian,english]{babel}
\usepackage{amsmath, amsfonts, amssymb, amsthm}

\usepackage[pdftex]{graphicx}
\graphicspath{{pictures/}}

\usepackage{pgf,tikz}
\usetikzlibrary{arrows}

\DeclareMathOperator{\normc}{g}
\DeclareMathOperator{\const}{const}

\newtheorem{theorem}{Theorem}
\newtheorem{proposition}{Proposition}

\theoremstyle{remark}

\theoremstyle{definition}
\newtheorem{definition}{Definition}
\newtheorem*{ack}{Acknowledgments}

\numberwithin{equation}{section}

\usepackage{caption}
\DeclareCaptionLabelSeparator{dot}{. }
\title[Isoperimetric property of $\lambda$-convex lunes]{About an isoperimetric property of $\lambda$-convex lunes on the Lobachevsky plane}
\author{Kostiantyn Drach}

\address{Geometry Department \\ V.N. Karazin Kharkiv National University\\ Svobody Sq. 4, 61022, Kharkiv \\ Ukraine}

\address{Department of Mathematical Analysis and Optimization \\ Sumy State University \\Rimskogo~- Korsakova str. 2, 40007, Sumy \\ Ukraine}
	
\email{drach@karazin.ua, kostya.drach@gmail.com}



\begin{document}

\maketitle

\begin{abstract}
We give a sharp lower bound on the area of a domain that can be enclosed by a closed embedded $\lambda$-convex curve of a given length on the Lobachevsky plane.
\end{abstract}

\section{Preliminaries and the main results}

The classical isoperimetric property of a circle in the two-dimensional space of constant curvature equal to $c$ claims that among all simple closed curves of a fixed length the maximal area is enclosed only by a circle. This property can be reformulated in an equivalent way in the form of an isoperimetric inequality. For an arbitrary simple closed curve of the length $L$ that encloses the domain of the area $A$ the following inequality holds (see, for example,~\cite{BZ})   
\begin{equation}
\label{isoineq}
L^2 - 4\pi A + c A^2 \geqslant 0,
\end{equation}
and the equality is attained only by circles.

Inequality~(\ref{isoineq}) gives a sharp upper bound on the area of the domain bounded by a curve provided that its length is fixed. At the same time, there exist simple closed curves that bound domains whose areas are arbitrary close to zero.  

A natural way of restriction the class of curves in order to obtain inequalities of different kind is to consider curves of bounded curvature. Such class appeared in a number of extremal problems (see, for example,~\cite{BorDr1, BorDr3, BorDr2, H, HTr, M}). In~\cite{BorDr1} for closed embedded $\lambda$-convex plane curves of a given length authors proved an inequality that gives a sharp lower bound on the area of domains enclosed by such curves. Similar inequality for curves on a sphere was obtained in~\cite{BorDr3}. In the present note we generalize these results for curves lying on a two-dimensional Lobachevsky plane $\mathbb H^2(-k^2)$ of Gaussian curvature equal to $-k^2$.

We recall the following definition.
\begin{definition} 
A locally convex curve $\gamma \subset \mathbb H^2(-k^2)$ is called \textit{$\lambda$-convex} with $\lambda \geqslant 0$ if for every point $P \in \gamma$ there exists a curve $\mu_P$ of constant geodesic curvature equal to $\lambda$ passing through $P$  in such a way, that in a neighborhood of $P$ the curve $\gamma$ lies from the convex side of $\mu_P$. 
\end{definition}

By the definition above, $0$-convex curves are just locally convex. It is known that on the Lobachevsky plane $\mathbb H^2(-k^2)$ there are three types of curves with constant geodesic curvature equal to $\lambda > 0$, namely, a circle (for $\lambda > k$), a horocycle (for $\lambda= k$), and an equidistant (for $k > \lambda > 0$).

We should also note that at  $C^r$-regular points of $\gamma$ with $r\geqslant 2$ the condition of being $\lambda$-convex is equivalent to the condition that at such points the geodesic curvature $\kappa_{\normc}$ of $\gamma$ satisfies the inequality $\kappa_{\normc} \geqslant \lambda$. Hence the class of $\lambda$-convex curves is a non-regular extension for the class containing smooth curves of geodesic curvature bounded from below by $\lambda$. 

It is known that a convex curve is twice continuously differentiable almost everywhere, and thus its geodesic curvature is almost everywhere well-defined. Therefore, a convex curve is $\lambda$-convex if and only if the inequality $\kappa_{\normc} \geqslant \lambda$ is satisfied at all points where geodesic curvature is defined.

For $\lambda > 0$ a \textit{$\lambda$-convex polygon} is a closed embedded $\lambda$-convex curve composed of arcs of curves with geodesic curvature equal to $\lambda$. It is known that there may be no-more-than countable number of such arcs.

A $\lambda$-convex polygon composed of two arcs of curves with curvature equal to $\lambda$ we will call a \textit{$\lambda$-convex lune} or simply a \textit{lune}.

It appears that the following theorem holds. 

\begin{theorem}
\label{mainth}

Let $\gamma$ be a closed embedded $\lambda$-convex curve (with $\lambda > 0$) lying on a two-dimensional Lobachevsky plane $\mathbb H^2 (-k^2)$ of Gaussian curvature equal to $-k^2$. If $L(\gamma)$ is the length of $\gamma$ and $A(\gamma)$ is the area of the domain enclosed by $\gamma$, then

\begin{enumerate}
\item
for $\lambda > k$ we have
\begin{equation}
\label{maintheq2}
A(\gamma) \geqslant \frac{\lambda}{k^2} L(\gamma) - \frac{4}{k^2} \arctan\left(\frac{\lambda}{\sqrt{\lambda^2 - k^2}}\tan\left(\frac{\sqrt{\lambda^2 - k^2}}{4}L(\gamma)\right)\right);
\end{equation}

\item
for $\lambda \geqslant k$ we have
\begin{equation}
\label{maintheq3}
A(\gamma) \geqslant \frac{1}{k} L(\gamma) - \frac{4}{k^2} \arctan\left(\frac{k}{4}L(\gamma)\right);
\end{equation}

\item
for $k > \lambda > 0$ we have
\begin{equation}
\label{maintheq4}
A(\gamma) \geqslant \frac{\lambda}{k^2} L(\gamma) - \frac{4}{k^2} \arctan\left(\frac{\lambda}{\sqrt{k^2 - \lambda^2}}\tanh \left(\frac{\sqrt{k^2 - \lambda^2}}{4}L(\gamma)\right)\right).
\end{equation}
\end{enumerate}
Moreover, the equality case in~(\ref{maintheq2})~-- (\ref{maintheq4}) holds only for $\lambda$-convex lunes. 
\end{theorem}

It is important to note that the inequality in Theorem~\ref{mainth} express the \emph{isoperimetric property} of $\lambda$-convex lunes. To make the statement above precise, we need to reformulate Theorem~\ref{mainth} in the following equivalent way.

\begin{theorem}
\label{mainthalt}
Let $\gamma$ be a closed embedded $\lambda$-convex curve (with $\lambda > 0$) lying on a two-dimensional Lobachevsky plane $\mathbb H^2 (-k^2)$. If $\gamma_0 \subset \mathbb H^2(-k^2)$ is a $\lambda$-convex lune such that $$L(\gamma) = L(\gamma_0),$$
then $$A(\gamma) \geqslant A(\gamma_0)$$ and the equality case holds if and only if $\gamma$ and $\gamma_0$ are congruent.
\end{theorem}

We will prove the main result in the form of Theorem~\ref{mainthalt}, and after that will show its equivalence to Theorem~\ref{mainth}.

\section{Proofs of the main results}

The principal tool for proving the main result is Pontryagin's Maximum Principle. We will follow a general approach from~\cite[\S 1.4]{MDmO}. 

In order to use Pontryagin's Maximum Principle we need to construct a controlled system. For this purpose let us introduce a so-called support function.

Let $O \in \mathbb H^2(-k^2)$ be a point inside a convex domain bounded by the curve $\gamma$. Let us consider on $\mathbb H^2(-k^2)$ the polar coordinate system with the origin at the point $O$ and with the angular parameter $\theta$ with $\theta \in [0, 2\pi)$. For each geodesic ray $OL$, emanating from $O$ and forming an angle $\theta$ with some fixed direction, let us consider a geodesic perpendicular to $OL$ and that is supporting for the curve $\gamma$ at some point $P$. By convexity of $\gamma$ such geodesic always exists, and the point $P$ is unique. Denote $h(\theta)$ to be the distance from the point $O$ to the supporting geodesic above, measured along the ray $OL$. The function $h(\theta) \colon [0, 2\pi) \to [0, +\infty)$ is a \textit{support function} for the curve $\gamma$. For our curve we will call a \textit{contact radius of curvature} the following quantity 
\begin{equation*}
g(\theta)=\frac{1}{k}\tanh\left(k h(\theta)\right).
\end{equation*}

We should note here that a strictly convex curve is uniquely determined by its support function. Remark also that the functions $h(\theta)$ and $g(\theta)$ belong to $C^{1,1}[0,2\pi]$ class of regularity. The last implies that $g(\theta)$ have a second derivative with respect to $\theta$ almost everywhere.

By a direct computation, which is similar to~\cite{Fil}, it can be easily shown that the contact radius of curvature $g$ is connected to the \textit{radius of curvature} $R(\theta)$, defined as ${1}/{\kappa_{\normc}(\theta)}$, by the following relation: 
\begin{equation}
\label{propsupfunceq1}
R = \frac{g'' + g}{\left(1-\frac{ k^2 {g'}^2}{1 - k^2 g^2}\right)^\frac{3}{2}} \text{ for almost all $\theta \in [0,2\pi]$,}
\end{equation}
(here by the prime sign we denote a derivative with respect the variable $\theta$).

In order to prove Theorem~\ref{mainthalt} let us fix the lengths of our curves and look for a minimum of the area of convex domains enclosed by these curves. To formalize this problem we need the expression for the length $L(\gamma)$ of a curve $\gamma$, and the expression for the area $A(\gamma)$ of the convex domain enclosed by $\gamma$ in terms of its contact radius of curvature $g(\theta)$ and its radius of curvature $R(\theta)$. Without loss of generality we may assume that the Lobachevsky plane has the Gaussian curvature equal to $-1$. Direct computations show that
\begin{equation}
\label{mainproofeq1}
L(\gamma) =  \int \limits_0^{2\pi} R \frac{\sqrt{1 - g^2 - {g'}^2}}{1 - g^2} d\theta, \,\, A(\gamma) = \int \limits_0^{2\pi} \left(\frac{\sqrt{1 - g^2 - {g'}^2}}{1 - g^2} - 1\right) d\theta.
\end{equation}

We remark here that the formulas above incorporate all possible jump angles at non-smooth points of $\gamma$.

Thus, we need to minimize $A(\gamma)$ taking into account~(\ref{propsupfunceq1}) and setting $L(\gamma) = \const$.
Let us interpret this problem as an optimal control problem with $t = \theta$ being a time variable, $u(t) = R(t)$ being a control, $x_1(t)=g(t)$, and $x_2(t) = \dot x_1(t) = g' (\theta)$ being phase variables.

Since $\gamma$ is a $\lambda$-convex curve, we have the restriction
\begin{equation}
\label{mainproofeq3}
0 \leqslant u(t) \leqslant \frac{1}{\lambda} \text{ a.e. on }[0, 2 \pi].
\end{equation}

Taking into consideration~(\ref{mainproofeq3}), and rewriting~(\ref{propsupfunceq1}), (\ref{mainproofeq1}) using the notations introduced above, we come to the following formal problem: 
\begin{equation}
\label{optcontrsys}
\begin{aligned}
&\int \limits_0^{2\pi} \left(\frac{\sqrt{1 - x_1^2 - x_2^2}}{1 - x_1^2} - 1\right) dt \rightarrow \min \\
&\int \limits_0^{2\pi} u \frac{\sqrt{1 - x_1^2 - x_2^2}}{1 - x_1^2} dt = \const \\
&\left\{
\begin{aligned}
&\dot x_1 = x_2 \\
&\dot x_2 = u \left(\frac{1 - x_1^2 - x_2^2}{1 - x_1^2}\right)^\frac{3}{2} - x_1  
\end{aligned} \right.\text{ a.e. on }[0, 2\pi]\\
&0 \leqslant u(t) \leqslant \frac{1}{\lambda} \text{\, a.e. on }[0, 2\pi]\\
&x_1(0) = x_1 (2 \pi) \\ &x_2(0) = x_2 (2 \pi).
\end{aligned}
\end{equation}

Moreover, in problem~(\ref{optcontrsys}) the control $u(t)$ is bounded measurable function on $[0,2 \pi]$, and the phase variable $x(t)$, defined as $x(t) = (x_1(t), x_2(t))$, is absolutely continuous function on $[0,2\pi]$ since $g(\theta) \in C^{1,1}[0,2\pi]$. In addition, all the functions used in the functional, the integral constraint and the controlled system are continuous with respect to all variables. The same smoothness condition holds for derivatives with respect to $x$ of the mentioned functions. 

Therefore, the pair $(x,u)$ that satisfies the controlled system from~(\ref{optcontrsys}) is a controlled process (see~\cite{MDmO}), and if $(x,u)$ also satisfies  the integral constraint and the boundary conditions of problem~(\ref{optcontrsys}), then the corresponding trajectory $\{(x(t),u(t)) \colon t\in[0,2\pi]\}$ is an admissible trajectory.

By Blaschke's selection theorem (see~\cite{Bla}) the posed problem of minimizing the area bounded by $\lambda$-convex curves while keeping their lengths fixed has a solution in the same class of curves. Hence the formalized version~(\ref{optcontrsys}) of the problem also has a solution. Thus in our case Pontryagin's Maximum Principle is a criterion for optimality of admissible trajectories. 

The adjoint system (see~\cite{MDmO}) for problem~(\ref{optcontrsys}) has the form
\begin{equation}
\label{adjp1}
\dot p_1 = p_2 \frac{(1 - x_1)^{\frac{5}{2}} + 3 u x_1 x_2^2 \sqrt{1 - x_1^2 - x_2^2}}{\sqrt{1 - x_1^2}} + \frac{x_1 (\mu_0 - \mu_1 u) \left(1 - x_1^2 - 2 x_2^2\right)}{\left(1 - x_1^2\right)^2 \sqrt{1 - x_1^2 - x_2^2}},
\end{equation}
\begin{equation}
\label{adjp2}
\dot p_2 = p_1 + p_2 \frac{3 u x_2 \sqrt{1 - x_1^2 - x_2^2}}{\left(1 - x_1^2\right)^{\frac{3}{2}}} - \frac{x_2 (\mu_0 - \mu_1 u)}{\left(1 - x_1^2\right) \sqrt{1 - x_1^2 - x_2^2}},
\end{equation}
where the adjoint variables $p_1(t)$ and $p_2(t)$ are absolutely continuous function on $[0,2\pi]$.

Also, Pontryagin's function for~(\ref{optcontrsys}) is equal to
\begin{equation}
\label{ham}
\begin{aligned}
\mathcal H (x, u, p, \mu_0, \mu_1) &= p_1 x_2 + p_2 \left(  u \left(\frac{1 - x_1^2 - x_2^2}{1 - x_1^2}\right)^\frac{3}{2} - x_1\right) \\
&+ \mu_1 \left(u \frac{\sqrt{1 - x_1^2 - x_2^2}}{1 - x_1^2}\right) - \mu_0 \left(\frac{\sqrt{1 - x_1^2 - x_2^2}}{1 - x_1^2} - 1\right),
\end{aligned}
\end{equation}
where $\mu_0$ and $\mu_1$ are some real numbers, and $\mu_0 \geqslant 0$. The variables $\mu_0$, $\mu_1$, $p_1$, and $p_2$ must satisfy the non-triviality condition (see~\cite{MDmO}).

Observe that Pontryagin's function~(\ref{ham}) is linear with respect to $u$, and can be written as $\mathcal H = u \mathcal H_1 + \mathcal H_2,$ where 
\begin{equation}
\label{hamh1}
\mathcal H_1 = \mu_1 \frac{\sqrt{1 - x_1^2 - x_2^2}}{1 - x_1^2} + p_2 \left( \frac{1 - x_1^2 - x_2^2}{1 - x_1^2}\right)^\frac{3}{2}.  
\end{equation}

From the maximality condition for $\mathcal H$ it follows that the optimal control for problem~(\ref{optcontrsys}) must be of the form 
\begin{equation}
\label{contr1}
u(t) = \left\{
\begin{aligned}
&{1}/{\lambda}, \text{ for } \mathcal H_1 > 0 \\
&0, \text{ for } \mathcal H_1 < 0 \\
&\text{undefined}, \text{ for } \mathcal H_1 = 0
\end{aligned}
\right. \text{ a.e. on }[0, 2\pi].
\end{equation}

In order to determine the control completely, we need to consider a so-called \textit{singular trajectory}, that is an admissible for~(\ref{optcontrsys}) trajectory on which $\mathcal H_1$ is identically zero for some interval $(t_1,t_2) \subset [0,2\pi]$.

From~(\ref{hamh1}) the condition $\mathcal H_1 = 0$ is equivalent to that on the singular trajectory
\begin{equation}
\label{singularextrp2}
p_2 = - \mu_1 \frac{\sqrt{1 - x_1^2} }{1 - x_1^2 - x_2^2}.
\end{equation}

Substituting~(\ref{singularextrp2}) into the differential equation~(\ref{adjp2}), and solving for $p_1(t)$ after all necessary cancellations, we will get that on the singular trajectory
\begin{equation}
\label{singularextrp1}
p_1 =  - x_2  \frac{\mu_1 x_1 \sqrt{1 - x_1^2} + \mu_0 \sqrt{1 - x_1^2 - x_2^2} }{(1 - x_1^2)(1 - x_1^2 - x_2^2)}.
\end{equation}

At the same time, the functions $p_1(t)$ and $p_2(t)$ have to satisfy the remaining equation~(\ref{adjp1}) from the adjacent system. If we substitute $p_1$ and $p_2$ into~(\ref{adjp1}) and simplify it using the expressions for $\dot x_1$ and $\dot x_2$ from~(\ref{optcontrsys}), then we will come to the equality 
\begin{equation*}
\label{adjsextrcons}
\frac{\mu_1 - \mu_0 u}{\left(1 - x_1^2\right)^{\frac{3}{2}}}=0.
\end{equation*} 
Hence on the whole interval $(t_1,t_2)$ the equality 
\begin{equation}
\label{conscond}
\mu_1 - \mu_0 u(t) \equiv 0
\end{equation}
must hold.

We can consider only so-called \textit{normal trajectroies} for which $\mu_0 = 1$. With this in mind, equality~(\ref{conscond}) is possible only if $u(t) \equiv \mu_1$ on the whole interval $(t_1,t_2)$. Since $\mu_1$ is a constant real number, and thus doesn't depend on the interval $(t_1,t_2)$, we conclude that the optimal control is equal $\mu_1$ on any arc of the singular trajectory.

Therefore,~(\ref{contr1}) can be rewritten as
\begin{equation*}
u(t) = \left\{
\begin{aligned}
1/\lambda,\,\, &\text{ for } \mathcal H_1 > 0 \\
0,\,\, &\text{ for } \mathcal H_1 < 0 \\
\mu_1, &\text{ for } \mathcal H_1 = 0
\end{aligned}
\right.
\end{equation*}

Let us show now that, in fact, in problem~(\ref{optcontrsys}) no arc of the singular extremal can be optimal. For such purpose let us consider the necessary condition for an arc of a singular trajectory to be a part of an optimal trajectory.

Recall that a natural number $q$ is an \textit{order} of the singular trajectory $(x^*,u^*)$ if it is a minimal number such that
$$\frac{\partial}{\partial u}\frac{d^{2q}}{dt^{2q}}\mathcal H_1 \left|_{(x^*,u^*)}\right. \neq 0, $$ 
where the time derivatives are taken with respect to the corresponding controlled system (see~\cite{KKM}).

The necessary condition for an arc of a singular trajectory of order $q$ to be a part of an optimal trajectory is so-called \textit{generalized Legendre~-- Clebsch condition} (see~\cite{KKM}, and also~\cite{B}), namely, along this arc
\begin{equation}
\label{LCC}
(-1)^q \frac{\partial}{\partial u} \left[\frac{d^{2q}}{dt^{2q}}\mathcal H_1\right] \leqslant 0.
\end{equation}

In our case we have the singular trajectory of order 1 ($q=1$). Indeed, using~(\ref{singularextrp2}) and~(\ref{singularextrp1}), along the singular trajectory we have
$$
-\frac{\partial}{\partial u}\frac{d^{2}\mathcal H_1}{dt^{2}} = \frac{\left(1 - x_1^2 - x_2^2\right)^\frac{3}{2}}{\left(1 - x_1^2\right)^3} > 0. 
$$

The last inequality contradicts the necessary  Legendre~-- Clebsch condition~(\ref{LCC}). Hence inclusion of the singular trajectory in the solution of problem~(\ref{optcontrsys}) is not optimal. Therefore, on the segment $[0, 2\pi]$ we have a bang-bang control with only values 0 or $1/\lambda$.

From the geometric point of view we obtained that an optimal curve must consist of arcs of curves with curvature equal to $\lambda$. Thus the solution of~(\ref{optcontrsys}) belongs to the class of $\lambda$-convex polygons with possibly infinite number of vertexes, that is points on a curve at which left and right semi-tangents do not coincide. For such a class we have the following geometric proposition, which proof is similar to the spherical case (see~\cite{BorDr3}, and also~\cite{PR}). 

\begin{proposition}
\label{mainthspecialcase}

Let $\gamma$ and $\tilde\gamma$ be, respectively, a $\lambda$-convex polygon and a $\lambda$-convex lune on a Lobachevsky plane $\mathbb H^2(-k^2)$. If $$L(\gamma) = L(\tilde\gamma),$$
then $$A(\gamma) \geqslant A(\tilde\gamma).$$
Moreover, the equality holds if and only if $\gamma$ and $\tilde\gamma$ are congruent.

\end{proposition}

Theorem~\ref{mainthalt} now follows directly from Proposition~\ref{mainthspecialcase}. At the same time, Theorem~\ref{mainth} follows from Theorem~\ref{mainthalt} and the relation between the length of a $\lambda$-convex lune and the area of the domain enclosed by it.

\begin{ack}
This work was partially done while the author was visiting the Centre de Recerca Matem\`atica as a participant of the Conformal Geometry and Geometric PDE's Program supported by a grant of the Clay Mathematics Institute. He would like to acknowledge both institutions for the given opportunities. The author is also thankful to prof. Alexander Borisenko for many fruitful discussions.
\end{ack}

\end{document}